\documentclass[11pt,leqno]{article} 
\usepackage{amsmath,amstext, amsthm, amssymb}

 \numberwithin{equation}{section}
\newtheorem{thm}{Theorem}[section]

\newtheorem{rem}[thm]{Remark}

\newcommand{\al}{\alpha}

\newcommand{\be}{\begin{equation}}
\newcommand{\ee}{\end{equation}}
\newcommand{\bea}{\begin{eqnarray}}
\newcommand{\eea}{\end{eqnarray}}

\newcommand{\ba}{\begin{array}}
\newcommand{\ea}{\end{array}}

\newcommand{\bt}{\beta}

\begin{document}

\title{Hypergeometric Origins of Diophantine Properties Associated With the Askey Scheme}
\author{Yang Chen\\
         Department of Mathematics\\
          Imperial College\\
          180 Queen's Gate\\
          London SW7 2BZ, UK\\
          ychen@ic.ac.uk\\
          \and Mourad E.H. Ismail
\\ Department of Mathematics \\ University of Central Florida
\\  Orlando, FL 32816\\
USA\\
ismail@math.ucf.edu}

\maketitle

\begin{abstract}
The "Diophantine" property of the zeros of certain polynomials in the Askey scheme, recently
discovered by Calogero and his collaborators, is explained, with suitably chosen parameter values,
in terms of the summation theorem of hypergeometric series. Here the Diophantine property refers to
integer valued zeros. It turns out that the same procedure can also
be applied to polynomials arising from the basic hypergeometric series. We found, with suitably chosen
parameters and certain $q-$analogue of the summation theorems, zeros of these polynomials explicitly, which are
no longer integer valued. This goes beyond the results obtained by the Authors mentioned above.
\end{abstract}


 \setcounter{equation}{0}
 \setcounter{thm}{0}
\section{Introduction}

 In a series of papers Calogero and his collaborators, see for example \cite{Bru:Cal:Dro1},
 \cite{Bru:Cal:Dro2}, and \cite{Bru:Cal:Dro3}, investigated various
integrable lattices of the Toda-type, with suitable boundary conditions. These
lattices arose as the dressing chains of Adler, Shabat, Yamilov and other. See for example
\cite{AdSh},\cite{ShYa} and \cite{Ya}.
It is found that if the
small amplitude motion about the equilibrium configuration is assumed to be isochronous, namely,
each component is periodic with the same period, then the characteristic frequencies, must
necessarily have integer values. Furthermore, if the assumption of nearest neighbor interaction is made
in the lattice models, then the secular equation whose zeros gives the characteristic frequencies reads
$$
\det(x\;\textsf{I}_N-\textsf{A}_N)=0,
$$
where $\textsf{A}_N$ is a tri-diagonal matrix of size  $N.$ We may take $N$ to be the number
of particles in the many-body problem. See \cite{Cal} for a detailed treatment. Hence
$$
P_N(x):=\det(x \textsf{I}_N-\textsf{A}_N)
$$
maybe interpreted as orthogonal polynomials if  the super-diagonal elements of
$\textsf{A}_N$ are real and none of them vanishes. We show that the Diophantine property are the generated
when the parameters of
the orthogonal polynomials are suitably chosen. The factorization occurs when the polynomials, although are
still characteristic polynomials of tri-diagonal matrices, are no longer orthogonal.

The motivation of considering this problem came from reading  \cite{Bru:Cal:Dro},
where a hypergeometric polynomial of degree $n$ is factored as $f_m(x)g_{n-m}(x)$, here $f_m$
has degree $m$, $g_{n-m}$ has degree $n-m$ and the zeros of $f_m$ are equi-spaced.
This holds for all $m, 1\le m \le n$. This factorization is referred to as
having the ``Diophantine property" in \cite{Bru:Cal:Dro}. Our explanation is that
all the Diophantine results in
\cite{Bru:Cal:Dro}  follow from summation theorems for hypergeometric functions and give $q$-analogues
 of all of them. This will be shown in \S 3. In \S 4 we provide $q$-analogues of all the results of
 \S 3, that is all the results in \cite{Bru:Cal:Dro}.
Section 2 contains the notation, summation theorems, and transformation formulas used in \S 3
and \S4.

 It is known that a sequence of monic orthogonal polynomials satisfy a three term recurrence relation
\bea
\label{eq3trr}
xP_n(x) = P_{n+1}(x) + \al_n P_n(x) + \bt_n P_{n-1}(x), n >0,
\eea
with $P_0(x):=1, P_1(x) :=x-\al_0$ then $P_n(x)$ can be represented as a determinant.
The monic polynomials have the determinant representation
\begin{equation}
\label{eq2.2.8}
P_n(x)=\begin{vmatrix}
x-\alpha_0  & -a_1       & 0          & \dotsm    &  0        & 0          & 0 \\
-a_1    & x-\alpha_1   & -a_2       & \dotsm    &  0        & 0          & 0 \\
\vdots  & \ddots  & \ddots  & \ddots    & \vdots    & \vdots     & \vdots \\
0        & 0        & \cdots     & {}        & -a_{n-2} & x-\alpha_{n-2} & -a_{n-1} \\
0        & 0        & 0          & \dotsm    &  0        & -a_{n-1}  & x-\alpha_{n-1}
\end{vmatrix},
\end{equation}
 where $a_n^2 = \bt_n, n >0$.   It is clear that if $\bt_k =0$ for some $k <n$ then $P_n(x)$ factors into a
 product of two polynomials of the same type, that is a product of two characteristic polynomials of
 tri-diagonal matrices. All the factorizations in the work of Bruschi, Calogero, and Droghei are of this type.
 What is surprising is that one of the two characteristic polynomials has equi-spaced zeros.

 The zeros of the polynomial $f_m$ in the  factorization of the Askey-Wilson polynomials
  turned out to be the points $[a q^k + q^{-k}/a]/2, k=0, 1, 2, \cdots$, where $a$ is one of
   the parameters in the Askey-Wilson polynomial. Such points, after $a \to ia$ are interpolation points
   in the sense that the values of an entire function  $h$ at these points determine the function uniquely
   provided that, $M(h,r)$, the maximum modulus of $h$ satisfies
   $M(h,r) \le Cr^\al \exp(b (\ln r)^2)$ with $- 2b\ln q  < 1$, for some $\al$, see \cite{Ism:Sta}.
   The integers are interpolation points for entire functions $h$ for which $M(h,r) \le C \exp(br)$ with $b < \pi$.

\begin{rem}
It is important to note that the Wilson polynomials are believed to be the most general orthogonal polynomials
of hypergeometric type  while the Askey-Wilson polynomials are
the most general orthogonal polynomials of basic hypergeometric type. As such we believe that it is
unlikely to extend this work to more general polynomials.
\end{rem}

 \section{Summation Theorems}
 Recall that the $q$-shifted factorial is
 \bea
 (a;q)_n = \prod_{k=1}^n(1-aq^{k-1}),
 \eea
 and a basic hypergeometric function is
   \bea
 \label{eqqHyps}
  \begin{gathered}
    {}_{r+1}\phi_r\left(\left. \begin{matrix}
a_1, a_2, \cdots, a_{r+1} \\
b_1, b_2, \cdots, b_r
\end{matrix}\, \right|q,z\right)
     = \sum_{n=0}^\infty \prod_{k=1}^{r+1} \frac{(a_k;q)_n}{( b_{k-1}; q)_n}\; z^n,
\end{gathered}
\eea
where $b_0:=q$.

 The Pfaff-Saalsch\"{u}tz theorem is,
  \bea
 \label{eqPfSa}
  \begin{gathered}
    {}_{3}F_2\left(\left. \begin{matrix}
-n, \quad  A, \quad B \\
C, 1+A+B -n - C
\end{matrix}\, \right|1\right)
     = \frac{(C-A)_n (C-B)_n}{(C)_n(C-A-B)_n}.
\end{gathered}
\eea
 The summation formula
 \bea
 \label{eqFieWimsu}
 \begin{gathered}
    {}_{r+2}F_{r+1}\left(\left. \begin{matrix}
 A, B, B_1+m_1, \cdots, B_r+m_r \\
B+1, B_1, \cdots, B_r
\end{matrix}\, \right|1\right) \\
= \frac{\Gamma(B+1)\Gamma(1-A)}{\Gamma(1+B-A)} \, \prod_{j=1}^r \frac{(B_j-B)_{m_j}}{(B_j)_{m_j}},
\end{gathered}
\eea
is known as the Karlsson-Minton sum, \cite{Gas:Rah}  but it follows from the earlier work of
 Fields and Wimp  \cite{Fie:Wim}. In particular we have
\bea
\label{eqFie-Wim}
 \begin{gathered}
    {}_{r+1}F_{r}\left(\left. \begin{matrix}
 A,  B_1+m_1, \cdots, B_r+m_r \\
 B_1, \cdots, B_r
\end{matrix}\, \right|1\right)
= 0,
 \end{gathered}
\eea
for Re $(-a) > m_1+m_2+ \cdots + m_r$.
We will also apply the Whipple transformation \cite{Sla}
\bea
\label{eqWhipple}
 \begin{gathered}
    {}_{4}F_{3}\left(\left. \begin{matrix}
 -n,   A, B, C \\
D, \; E,\;    F
\end{matrix}\, \right|1\right)
= \frac{(E-A)_n(F-A)_n}{(E)_n(F)_n} \\
\times  {}_{4}F_{3}\left(\left. \begin{matrix}
 -n, \;\;A, \;\;  D-B, \; \;  D- C \\
D, A+1-n-E, A+1-n-F
\end{matrix}\, \right|1\right)
 \end{gathered}
\eea
where $D+E+F= A+B+C+1-n$.

The $q$-analogue of the Pfaff-Saalsch\"{u}tz theorem is
\bea
 \label{eqq-PfSa}
  \begin{gathered}
    {}_{3}\phi_2\left(\left. \begin{matrix}
q^{-n}, \quad  A, \quad B \\
C, q^{1-n}AB/C
\end{matrix}\, \right|q, q \right)
     = \frac{(C/A;q)_n (C/B;q)_n}{(C;q)_n(C/AB;q)_n},
\end{gathered}
\eea
while the $q$-analogue of the Whipple transformation is the Sears transformation \cite[(III.15)]{Gas:Rah}
is
\bea
 \begin{gathered}
    {}_{4}\phi_{3}\left(\left. \begin{matrix}
 q^{-n},   A, B, C \\
D, \; E,\;    F
\end{matrix}\, \right|q,q\right)
=  A^n\, \frac{(E/A;q)_n(F/A;q)_n}{(E;q)_n(F;q)_n} \\
\qquad\qquad  \times  {}_{4}\phi_{3}\left(\left. \begin{matrix}
 q^{-n}, \;\;A, \;\;  D/B, \; \;  D/ C \\
D, q^{1-n}A/E, q^{1-n}A/F
\end{matrix}\, \right|q,q\right)
 \end{gathered}
 \label{eqSears}
\eea
where $DEF= q^{1-n}ABC$.

Some useful identities are:
\begin{gather}
(aq^{-n};q)_n = (q/a;q)_n \left(-a\right)^n\, q^{-\binom{n+1}{2}},
\label{eqqshift1}\\
(aq^{-n};q)_{n-k}= \frac{(q/a;q)_n}{(q/a;q)_k}\left(-a\right)^{n-k}\,
q^{\binom{k+1}{2}-\binom{n+1}{2}}.  \label{eqqshift2}
\end{gather}
Of course the first is a special case of the second.

 \section{Complete factorization of the Wilson and Related Polynomials}

 The Wilson polynomials is
 \bea
 \begin{gathered}
 W_n(x;{\bf t}) = \prod_{j=1}^3(t_1+ t_j)_n  \qquad \qquad \\
\qquad  \times  {}_{4}F_{3}\left(\left. \begin{matrix}
 -n, t_1+t_2+t_3+t_4 + n-1, t_1+i\sqrt{x}, t_1-i\sqrt{x} \\
 t_1+t_2, t_1+t_3, t_1+t_4
\end{matrix}\, \right|1\right),
\end{gathered}
 \eea
where ${\bf t}:=(t_1,t_2, t_3, t_4)$. It is a fact that the Wilson polynomials is  symmetric in the
 parameters $t_1,t_2, t_3, t_4$. The invariance of $W_n$ under permutations of $\{t_2,t_3,t_4\}$
 is obvious but the invariance under permuting $t_1$ and $t_j$, for $j =2,3,4$ is not obvious and
 is called the Whipple transformation, \cite{Sla}.

If we wish to  find a complete factorization of $W_n$, or equivalently identify  all the zeros of  $W_n$,
then we must choose the parameters in such a way that the ${}_4F_3$ representation can be summed explicitly.

We shall denote the monic Wilson polynomials by $\{\tilde{W}_n(x;{\bf t})\}$, that is
\bea
\begin{gathered}
\tilde{W}_n(x;t_1,t_2,t_3,t_4)  \qquad \qquad
  \\= \frac{(-1)^n}{(n+t_1+t_2+t_3+t_4-1)_n} \, W_n(x;t_1,t_2,t_3,t_4).
\end{gathered}
\eea

\bigskip

\noindent{\bf Case 1}.  We reduce the ${}_4F_3$  to a ${}_3F_2$ and use the Pfaff-Saalsch\'{u}tz theorem.
 Since we want to keep $x$ in the factorization we demand that $n-1+\sum_{k=1}^4 t_k$ be equal to $t_1+t_j$
  for some $j$. It is easy to see that this happens if and only if $t_i+t_j = 1-n$ for some $i\ne j, 1 < i,j \le 4$.
   There is no loss of generality in assuming $t_4 = 1-n-t_3$. In this case \eqref{eqPfSa} gives
\bea
\begin{gathered}
\tilde{W}_n(x;t_1,t_2,t_3,1-n-  t_3)=(-1)^n(t_1+t_3)_n(t_1+1-n-t_3)_n\\
\times {}_{3}F_2\left(\left. \begin{matrix}
-n, t_1+i\sqrt{x},  t_1- i\sqrt{x} \\
t_1+t_3, t_1+1-n-t_3
\end{matrix}\, \right|1\right)\\
= (-1)^n(t_1+t_3)_n(t_1+1-n-t_3)_n \, \times\, \frac{(t_3+i\sqrt{x})_n(t_3-i\sqrt{x})_n}
{(t_1+t_3)_n(t_3-t_1)_n}.
\end{gathered}
\nonumber
\eea
Since $(-1)^n(t_1-t_3+1-n)_n = (t_3-t_1)_n$, it follows that
\bea
\begin{gathered}
\tilde{W}_n(x;t_1,t_2,t_3,1-n-  t_3)= (t_3+i\sqrt{x})_n(t_3-i\sqrt{x})_n  \\
= \prod_{k=1}^n[x^2+(t_3+k-1)^2].
\end{gathered}
\eea
The above factorization is  equations (33) and (41a) of \cite{Bru:Cal:Dro}.

\bigskip

\noindent{\bf Case 2}.  We identify values of  $x$ suitable to apply
\eqref{eqFie-Wim}. For example, we may choose $i\sqrt{x} = t_4 +j$, and
make $t_1-i\sqrt{x}$, which is $t_1 -t_4-j$, equal $t_1+ t_3+k$. Thus we make the parameter identification
$x = -(t_4+j)^2, t_3 = -t_4-j-k$. Finally we demand that $n + t_3+t_4-1 = s$, which is equivalent
to $n-1 \ge j+k$.  Now we set $m=j+k+1$, replace $j$ by $j-1$ and we find that
\bea
\label{eqcase2}
W_n(-(t_4+j-1)^2; t_1, t_2, -t_4-m +1, t_4) =0,
\eea
for $1\le j \le m, 1\le m \le n$. This is (35) in \cite{Bru:Cal:Dro}. This is particularly
interesting because it seems to give a partial factorization when $m <n$. We shall return to
this point at the end of this section.

When $m=n$ we obtain the factorization
\bea
\tilde{W}_n(x; t_1, t_2, -t_4-n +1, t_4)= \prod_{j=1}^n[x+ (t_4+ j-1)^2],
\eea
which is (34) in \cite{Bru:Cal:Dro}. The special case $t_4 = (1-2n)/4$ is (41a) of
 \cite{Bru:Cal:Dro}. This last factorization also follow from the Pfaff-Saalsch\'{u}tz formulas
 since the ${}_4F_3$ reduces to a ${}_3F_2$. Indeed in this case we have
 \bea
    \begin{gathered}
 W_n(x; t_1, t_2, -t_4-n +1, t_4) = (t_1+t_2)_n(t_1+t_4)_n(t_1+ 1-n-t_4)_n \\
\times     {}_{3}F_2\left(\left. \begin{matrix}
-n, t_1+i\sqrt{x}, t_1-i\sqrt{x} \\
 t_1+t_4, t_1+1-n-t_4
\end{matrix}\, \right|1\right)\\
     =   (t_1+t_2)_n(t_1+t_4)_n(t_1+ 1-n-t_4)_n\frac{(t_4+i\sqrt{x})_n(t_4-i\sqrt{x})_n}
     {(t_1+t_4)_n(t_4-t_1)_n}.
\end{gathered}
\nonumber
\eea

 We now show how to discover \eqref{eqcase2}
 from the Whipple transformation.  It is clear that $W_n(x; t_1, t_2, 1-t_4-m, t_4)$ is a
 constant multiple of
 \bea
 \begin{gathered}
 \label{eqfirst4F3}
 {}_{4}F_{3}\left(\left. \begin{matrix}
 -n, t_1+t_2+ n -m, t_1+i\sqrt{x}, t_1-i\sqrt{x} \\
 t_1+t_2, t_1+1-m-t_4, t_1+t_4
\end{matrix}\, \right|1\right) \\
= \frac{(1-m-t_4-i\sqrt{x})_n(t_4-i\sqrt{x})_n}{(t_1-t_4+1-m)_n(t_1+t_4)_n} \\
\quad \times   {}_{4}F_{3}\left(\left. \begin{matrix}
 -n, \; \; \;  t_1+i\sqrt{x}, \; \; \; m-n, \; \; \; t_2+ i\sqrt{x} \\
 t_1+t_2, t_4 +i\sqrt{x} +m-n, i\sqrt{x} +1 -n -t_4
\end{matrix}\, \right|1\right).
 \end{gathered}
\eea
In the last step we applied the Whipple transformation \eqref{eqWhipple} with the parameter identification
\bea
\nonumber
 \begin{gathered}
A = t_1+i\sqrt{x}, B = n +t_1+t_2 -m,  C = t_1- i\sqrt{x}, \\
D = t_1+t_2, E = t_1-t_4+1-m, F = t_1+t_4.
 \end{gathered}
\eea
We next apply the Whipple transformation again with the choices
\bea
\nonumber
 \begin{gathered}
A = t_2+i\sqrt{x}, B = t_1+i\sqrt{x},  C = -n, \\
D = t_1+t_2, E = t_4 + i \sqrt{x} +m-n, F = 1-n-t_4+i\sqrt{x},
 \end{gathered}
\eea
and $n$ is now $n-m$. Therefore  the left-hand side of \eqref{eqfirst4F3}
is
\bea
 \begin{gathered}
  \frac{(1-m-t_4-i\sqrt{x})_n(t_4-i\sqrt{x})_n}{(t_1-t_4+1-m)_n(t_1+t_4)_n} \\
\quad \times   \frac{(t_4-t_2+m-n)_{n-m}(1-n-t_2-t_4)_{n-m}}
   {(t_4+m-n+i\sqrt{x})_{n-m}(1-n-t_4+i\sqrt{x})_{n-m}}\\
\quad \times   {}_{4}F_{3}\left(\left. \begin{matrix}
 m-n, \; \; \; n+t_1+t_2, \;\;\; t_2+i\sqrt{x}, \; \; \; t_2- i\sqrt{x} \\
 t_1+t_2, t_2+1- t_4, t_2+t_4+m \end{matrix}\, \right|1\right).
 \end{gathered}
\eea
It is clear that the ${}_4F_3$ in the above expression  is a constant multiple of
$W_{n-m}(x; t_2, t_1, 1-t_4, m+t_4)$. We now apply  the identities
\bea
(\alpha+1)_n = (-1)^n(-\alpha-n)_n,  \quad (\alpha)_n = (\alpha)_m(\alpha+m)_{n-m}
\eea
to see that
\bea
\begin{gathered}
\frac{(1-m-t_4-i\sqrt{x})_n(t_4-i\sqrt{x})_n}
 {(t_4+m-n+i\sqrt{x})_{n-m}(1-n-t_4+i\sqrt{x})_{n-m}} \\
 = \frac{(1-m-t_4-i\sqrt{x})_n(t_4-i\sqrt{x})_n}
   {(1-t_4-i\sqrt{x})_{n-m}(m+t_4-\sqrt{x})_{n-m}}\\
   = \frac{(1-m-t_4-i\sqrt{x})_m(t_4-i\sqrt{x})_n}
   {(m+t_4-\sqrt{x})_{n-m}}\\
   = (-1)^m(t_4+i\sqrt{x})_m(t_4-i\sqrt{x})_m = (-1)^m\prod_{j=0}^{m-1}[x^2+(t_4+j)^2].
    \end{gathered}
    \nonumber
\eea
This shows that
\bea
\begin{gathered}
\frac{W_n(x; t_1,  t_2, 1-t_4-m , t_4)}{(t_1+t_2)_n(t_1+ t_2, 1-t_4-m)_n(t_1+t_4)_n}
\\
=\frac{ (-1)^m(t_4+i\sqrt{x})_m(t_4-i\sqrt{x})_m}
 {( t_1+t_2)_{n-m}(t_1+t_4)_{n-m}(t_1- t_4+ 1-m)_{n-m}}\\
 \quad \times
 W_{n-m}(x; t_2,  t_1, 1-t_4 , t_4+m).
    \end{gathered}
    \nonumber
\eea

\section{The Askey-Wilson Polynomials}
The Askey-Wilson polynomials are defined through the representation
\bea
 \begin{gathered}
 p_n(x;{\bf t}) = t_1^{-n}\; \prod_{j=1}^3(t_1 t_j;q)_n  \qquad \qquad \\
\qquad  \times  {}_{4}\phi_{3}\left(\left. \begin{matrix}
 q^{-n}, t_1t_2t_3t_4 q^{n-1},  t_1 e^{i\theta}, t_1e^{-i\theta} \\
 t_1t_2, t_1t_3, t_1t_4
\end{matrix}\, \right|q,q\right).
\end{gathered}
 \eea
As in Case 1 of \S 3 we reduce the ${}_4\phi_3$ to a ${}_3\phi_2$.
There is no loss of generality in assuming that
$q^{n-1}t_1t_2 t_3t_4 = t_1t_2$, that is $t_4= q^{n-1}/t_3$. Applying the summation formula
\eqref{eqq-PfSa} we get
\bea
\begin{gathered}
p_n(\cos \theta; t_1, t_2, t_3, q^{n-1}/t_3) = \\
(t_3e^{i\theta};q)_n(t_3e^{-i\theta};q)_n\frac{(t_1t_2;q)_n q^{1-n}t_1/t_3;q)_n}
{t_1^n(q/t_1t_3;q)_n}
\end{gathered}
\eea
 It is clear that
 \bea
 (t_3e^{i\theta};q)_n(t_3e^{-i\theta};q)_n = \prod_{k=1}^n[1-\, 2t_3xq^{k-1}\, +t_3^2\, q^{2k-2}],
 \eea
 from which we can find the zeros explicitly.

As in of  Case 2 \S 3  we choose  $t_3 = q^{1-m}/t_4$ then apply the Sears transformation
\eqref{eqSears} with the parameter identification
\bea
\nonumber
 \begin{gathered}
A = t_1 e^{i\theta} , \quad B = q^{n-m} t_1t_2,   \quad  C = t_1e^{-i\theta}, \\
D = t_1t_2,  \quad E =  q^{1-m}t_1/t_4,  \quad  F = t_1t_4.
 \end{gathered}
\eea
The result is
\bea
\nonumber
 \begin{gathered}
 p_n(\cos \theta;{\bf t}) = e^{in\theta} (t_1t_2;q)_n(q^{1-m}e^{i\theta}/t_4;q)_n(t_4e^{-i\theta};q)_n
\\
\qquad  \times  {}_{4}\phi_{3}\left(\left. \begin{matrix}
 q^{-n},  t_1 e^{i\theta},  q^{m-n},  t_2e^{i\theta} \\
 t_1t_2, q^{m-n}t_4e^{i\theta}, q^{1-n}e^{i\theta}/t_4
\end{matrix}\, \right|q,q\right).
\end{gathered}
 \eea
We next apply the Sears  transformation again with the choices
\bea
\nonumber
 \begin{gathered}
A = t_2e^{i\theta},  \quad  B = t_1 e^{i\theta},   \quad  C = q^{-n}, \\
D = t_1t_2,  \quad  E = t_4 q^{m-n} e^{i\theta},  \quad  F = q^{1-n}e^{i\theta}/t_4,
 \end{gathered}
\eea
and the terminating parameter $n$ is replaced by $n-m$. This leads to
\bea
\nonumber
 \begin{gathered}
 p_n(x;{\bf t})
 =\frac{t_2^{2(n-m)}(t_1t_2;q)_n(q^{1-m}e^{-i\theta}/t_4;q)_n
 (t_4e^{-i\theta};q)_n}{(t_1t_2;q)_{n-m}(t_2q/t_4;q)_{n-m}(q^m t_2t_4;q)_{n-m}}
\\
  \times \frac{(q^{m-n}t_4/t_2;q)_{n-m}(q^{1-n}/t_2t_4;q)_{n-m}}
  {(q^{m-n}t_4 q^{i\theta};q)_{n-m}(q^{1-n}e^{i\theta}/t_4;q)_{n-m}}
  e^{i(2n-m)\theta}\\
  \times
  p_{n-m}(x; t_2, t_1, q/t_4, q^m t_4),
\end{gathered}
 \eea
 with $x = \cos \theta$.
Applying equations \eqref{eqqshift1}--\eqref{eqqshift2} we finally establish the factorization
\bea
 \begin{gathered}
 p_n(x;{\bf t}) =  (t_4e^{i\theta};q)_m (t_4e^{-i\theta};q)_m \, p_{n-m}(x; t_2, t_1, q/t_4, q^m t_4) \\
 \times
(-1)^m \,t_4^{n-2m}t_2^{m-n} \frac{(q^m t_1t_2;q)_{n-m}}{(q^m t_1t_4;q)_{n-m}},
\end{gathered}
\eea
again with $x = \cos \theta$.
Note that
\bea
 (t_4e^{i\theta};q)_m (t_4e^{-i\theta};q)_m = \prod_{k=1}^n[1-2xt_4q^{k-1}+ t_4^2q^{2k-2}].
\eea

\bigskip


{\bf Acknowledgments}.  The research of Mourad Ismail is supported by a research
grant from King Saud University, Riyadh, Saudi Arabia. Part of this paper was written while
the second  author was visiting the Isaac Newton Institute as part of the discrete integrable
 systems  program. He wishes
to thank the Institute's staff and the organizers for the hospitality and the excellent scientific environment.

\end{document}